\documentclass[11pt]{article}
\usepackage{amsmath,amssymb,enumerate} 

\usepackage{amsthm}
\usepackage[subnum]{cases}

\usepackage{tikz}
\usetikzlibrary{arrows}
\usetikzlibrary{positioning}

\newtheorem{theorem}{Theorem}
\newtheorem{lemma}{Lemma}

\newtheorem{example}{Example}
\newtheorem{remark}{Remark}

\def\cH{\mathcal{H}}

\def\cF{\mathcal{F}}

\def\G{\mathcal{G}}

\def\P{\mathcal{P}}
\def\N{\mathcal{N}}

\def\U{\mathcal{U}}

\def\ZZP{\mathbb{Z}_{+}}
\def\ZZ{\mathbb{Z}}
\def\ZZ1{\mathbf{Z}_{\geq 1}}

\def\be{\mathbf{e}}
\def\ba{\mathbf{a}}
\def\bb{\mathbf{b}}
\def\bx{\mathbf{x}}
\def\bc{\mathbf{c}}
\def\by{\mathbf{y}}

\def\M{M}
\def\Y{Y}
\def\H{h_\cH}

\def\Nim{\mbox{{\sc Nim}}}

\title{On the Sprague-Grundy function of 
compound games}

\author{
Endre Boros\thanks{MSIS and RUTCOR, RBS, Rutgers University,
100 Rockafeller Road, Piscataway, NJ 08854; e-mail:
endre.boros@rutgers.edu}
\and
Vladimir Gurvich\thanks{
National Research University Higher School of Economics (HSE) Moscow Russia; e-mail:
vgurvich@hse.ru}
\and
Levi Kitrossky\thanks{
Mobileye, An Intel Company
13 Hartom St. Jerusalem 9777513 Israel; e-mail:
lkitrossky@gmail.com}
\and
Kazuhisa Makino\thanks{
Research Institute for Mathematical Sciences $($RIMS$)$ Kyoto University,  
Kyoto 606-8502, Japan; e-mail: makino@kurims.kyoto-u.ac.jp} 
}

\begin{document}

\maketitle

\abstract{
The classical game of  {\sc Nim}  can be naturally extended
and played on an arbitrary hypergraph  $\cH \subseteq 2^V \setminus \{\emptyset\}$
whose vertices $V = \{1, \ldots, n\}$  correspond to piles of stones.
By one move a player chooses an edge  $H$  of  $\cH$ and
reduces arbitrarily all piles  $i \in H$.
In 1901 Bouton solved the classical  {\sc Nim} 
for which  $\cH = \{\{1\}, \ldots, \{n\}\}$.
In 1910 Moore introduced and solved a more general game  $k$-{\sc Nim},
for which  $\cH = \{H \subseteq V \mid |H| \leq k\}$, where $1 \leq k < n$.
In 1980 Jenkyns and Mayberry obtained an explicit formula
for the Sprague-Grundy function of Moore's  {\sc Nim}  for the case  $k+1 = n$.
Recently it was shown that the same formula works for
a large class of hypergraphs.
In this paper we study combinatorial properties of these hypergraphs and   
obtain explicit formulas for the Sprague-Grundy functions
of the conjunctive and selective compounds of the corresponding hypergraph {\sc Nim}  games.
\newline

\newline 
MSC classes: 91A46, 91A05
}

\section{Introduction}
\label{s1}

In the classical game of {\sc Nim} there are
$n$  piles of stones and two players move alternately.
A move consists of choosing a nonempty pile and taking some positive number of stones from it.
The player who cannot move is the loser.
Bouton \cite{Bou901} analyzed this game and described the winning strategy for it.

In this paper we consider the following generalization of {\sc Nim}.
For a positive integer $n$, let us denote by $V = \{1,...,n\}$ a set of $n$ piles of stones. 
Let $\ZZP$ denote the set of nonnegative integers.
We use $x\in\ZZP^V$ to describe a position, where coordinate $x_i$ denotes
the number of stones in pile $i \in V$.
Given a hypergraph $\cH \subseteq 2^V$ a move from a position $x\in\ZZP^V$
consists in choosing an edge $H \in \cH$ and strictly decreasing all $x_i$ values for $i\in H$.
The game starts in an initial position $x\in\ZZP^V$ and involves two players who alternate in making moves.
Similarly to {\sc Nim}, the player who cannot move is the loser.
Such games were considered in \cite{BGHMM15,BGHMM16,BGHMM18} and called hypergraph {\sc Nim}.
We denote by {\sc Nim}$_\cH$ an instance of this family.
We assume in this paper that  
$V=\bigcup_{H \in \cH} H$  and $\emptyset \not\in \cH$
for all considered hypergraphs  $\cH\subseteq 2^V$.
In other words, every move strictly decreases some of the piles.

Hypergraph {\sc Nim} games are impartial.
In this paper we do not need to immerse in the theory of impartial games. We
recall only a few basic facts to explain and motivate our research.
We refer the reader to \cite{Alb07,BCG01-04,Sie13} for more details.

A position of an impartial game is called \textit{winning}, or an $\N$-\textit{position},
if starting from it
the first player can win, no matter what the second player does.
The remaining positions are called \textit{losing}, or $\P$\textit{-positions}.
It is known that every move from a $\P$-position goes to an $\N$-position,
while from any $\N$-position there always exists a move to a $\P$-position.
The so-called Sprague-Grundy (SG) function $\G_\Gamma$ of an impartial game $\Gamma$ is
a refinement of the above $\P$-$\N$ partition, see Section \ref{section-example} for the definition.
Namely, $\G_\Gamma(x)=0$ if and only if $x$ is a $\P$-position.
The notion of the SG function was introduced independently by Sprague and Grundy
\cite{Spr35, Spr37, Gru39} and it plays a fundamental role
in solving \emph{disjunctive sums} of impartial games.

Finding a formula for the SG function of an impartial game remains a challenge.
Closed form descriptions are known only for some special cases.
We recall below some known results.
The purpose of our research is to extend these results and
to describe classes of hypergraphs for which
we can provide a closed formula for the SG function of {\sc Nim}$_\cH$.

The game {\sc Nim}$_\cH$ is a common generalization
of several families of impartial games considered in the literature.
For a subset $S \subseteq V$ and a positive integer
$k \leq |S|$  we denote by
\[
\binom{S}{k}=\{H\subseteq S\mid |H|=k\}.
\]
For instance, if $\cH=\binom{V}{1}$ then {\sc Nim}$_\cH$ is the classical {\sc Nim}.
The case of $\cH=\bigcup_{j=1}^k \binom{V}{j}$, where $k<n$, was considered by Moore \cite{Moo910}.
He characterized for these games the set of $\P$-positions, that is those with SG value $0$.
Jenkyns and Mayberry \cite{JM80} described also the
set of positions in which the SG value is $1$ and
provided an explicit formula for the SG function in case of $k = n-1$.
This result was extended in \cite{BGHM15}.
In \cite{BGHMM15} the game {\sc Nim}$_\cH$ was considered
in the case of $\cH=\binom{V}{k}$
and the corresponding SG function was determined when  $2k\geq n$.
Further examples such as matroid,  $2$-uniform (graph), symmetric, and hereditary hypergraph {\sc Nim} games are 
considered in \cite{BGHMM16,BGHMM17,BGHMM18,ES96}.
Surprisingly, for many of these examples the SG function is described by the same formula,
a special case of which was introduced by Jenkyns and Mayberry \cite{JM80}.
In honor of their contribution, this formula and the hypergraphs
for which it describes the SG function were called \emph{JM} in \cite{BGHMM18}.

In this paper, we consider compositions of games and their SG functions for JM hypergraph {\sc Nim} games.  
There are three types of compounds considered in the literature.
Given two games  $\Gamma_1$ and  $\Gamma_2$
with disjoint sets  of positions $X_1$ and $X_2$,
the compound game  $\Gamma$   has the set of  positions
$X = X_1 \times X_2$, while  the set  of its moves
can be introduced in three different ways as follows.

\begin{description}
\item[Disjunctive compound  $\Gamma_1 \oplus \Gamma_2$:] 
a player makes a move in exactly one of the two games:
either in  $\Gamma_1$  or  in  $\Gamma_2$.
\item[Conjunctive compound $\Gamma_1 \otimes \Gamma_2$:] 
a player makes a move in both games: one in  $\Gamma_1$  and one in  $\Gamma_2$.
\item[Selective compound $\Gamma_1 \boxplus \Gamma_2$:] 
a player makes a move either in one of the two games or in both.
\end{description}

Obviously, all three operations $\oplus$, $\otimes$, and $\boxplus$ are
associative and commutative and, hence all three compounds are
well defined not  only for two, but for any number of component games.

\medskip

The disjunctive compound was introduced by Sprague and  Grundy  \cite{Spr35,Spr37,Gru39};
the conjunctive and selective ones were  added by Smith and Conway \cite{Smi66,Con76}.
In \cite{BGHMM16} a concept of hypergraph combination of games,
which generalizes all three above compounds,  was introduced.

To state our main results we need a few more definitions.
To integers $m$, $y$ and $h$ let us associate the following quantities:
\begin{equation}\label{e-vmy}
v(m,y) ~=~ \binom{y+1}{2} +
\left( \left(m-\binom{y+1}{2}-1\right)\mod \left(y+1\right)\right)
\end{equation}
and
\begin{numcases}{\U(m,y,h) =}
h  \quad \text{ if } m \leq \binom{y+1}{2}\label{e-long}\\
v(m,y) \quad \text{ otherwise.} \label{e-short}
\end{numcases}

Given a hypergraph $\cH\subseteq 2^V$, the {\em height} $h_\cH(\bx)$ of a position $\bx\in\ZZP^V$ 
is defined as the maximum number of consecutive moves that
the players can make in {\sc Nim}$_\cH$ starting from $\bx$.
Furthermore, for a position $\bx \in \ZZP^V$ of $\Nim_\cH$ we define
\begin{subequations}
\label{e-myzv}
\begin{align}
m(\bx) &=\min_{i\in V} x_i, ~~\text{ and } \label{e-m}\\
y^{}_{\cH}(\bx) &=h_\cH(\bx-m(\bx)\be) \label{e-y}
\end{align}
\end{subequations}
where $\be=(1,1,...,1)$ is the $n$-vector of ones.
A position $\bx$ is called {\em long} if $m(\bx)\leq \binom{y^{}_{\cH}(\bx)+1}{2}$ and it is called {\em short} otherwise.

The expression
$\U(\bx)=\U(m(\bx), y_\cH(\bx), h_\cH(\bx))$ for a position $\bx\in\ZZP^V$ is called the JM formula.  
We call the hypergraph $\cH$ JM if the JM formula  represents the SG function of {\sc Nim}$_\cH$. 

In this paper we focus on a special subclass of JM hypergraphs.

\begin{itemize}
\item[(A)] Given a hypergraph $ \cH\subseteq 2^V$,
an edge $H \in \cH$ is called a \emph{transversal edge}
if it intersects every edge of the hypergraph, that is, if
$H \cap H' \neq \emptyset$ for all $H' \in \cH$.
A hypergraph with no transversal edge is called \emph{transversal-free}.
For a subset $S\subseteq V$ we denote by $\cH_S=\{H\in\cH \mid H\subseteq S\}$ the induced subhypergraph.
A hypergraph $\cH$ is called \emph{minimal transversal-free}
if it is transversal-free, but any proper induced subhypergraph of it has a transversal edge.

\item[(B)] Let us call a hypergraph $\cH$ \emph{minimum-decreasing}
if for every position $\bx \in \ZZP^V$ of {\sc Nim}$_\cH$
there exists a  move $\bx \to \bx'$
such that $m(\bx')< m(\bx)$,  $h_{\cH}(\bx') = h_{\cH}(\bx)-1$, and $x_i-x'_i \leq 1$ for all $i \in V$.

\item[(C)] A sequence of edges $H_0, H_1,\dots ,H_q$ in $\cH$ is called a \emph{chain}, if
$H_{i+1}\cap H_i\neq \emptyset$,    $|H_{i+1}\setminus H_i|=1$, and $H_i \subseteq H_0 \cup H_q$  for all $i=0,1,\dots,q-1$.
We say that a hypergraph $\cH$ has the \emph{chain-property}
if for any two distinct edges $H, H' \in \cH$  
there exists a chain $H_0,\dots,H_q$ in $\cH$ such that $H = H_0$ and $H' = H_q$.
\end{itemize}

A hypergraph satisfying properties (A), (B), and (C)
will be called a JM+ hypergraph. It was shown in \cite{BGHMM17}
that every JM+ hypergraph is JM, and that property (A) is necessary for a hypergraph to be JM.
We show in Section \ref{s6} that no two of the above three properties imply that the hypergraph is JM,
and in particular they do not imply the third property. We also show that unlike (A), property (C) is not necessary for a hypergraph to be JM.
The necessity of (B) remains an open question.

Let us note that JM+ is a proper subfamily of JM, since for instance some of the symmetric JM hypergraphs described in \cite{BGHMM18} do not belong to JM+. On the other hand, JM+ contains all JM hypergraphs described in \cite{BGHMM17}, including JM matroids and JM graphs.
It seems to be a challenging open problem to find a combinatorial characterization of JM hypergraphs.

Let us consider hypergraphs $\cH_i\subseteq 2^{V_i}$, $i=1,...,p$, where the sets $V_i$, $i=1,...,p$ are pairwise disjoint,
and define
\begin{align}
\cH_1\otimes\dots\otimes\cH_p &=~ \left\{ \left.\bigcup_{i=1}^p H^i\right| H^i\in\cH_i,~ i=1,...,p \right\}, \text{ and }\\
\cH_1\boxplus\dots\boxplus\cH_p &=~ \left\{ \left.\bigcup_{i=1}^p H^i\right| H^i\in\cH_i\cup\{\emptyset\},~ i=1,...,p \right\} \setminus \{\emptyset\}.
\end{align}
We call $\cH_1\otimes\dots\otimes\cH_p$ the conjunctive compound, and $\cH_1\boxplus\dots\boxplus\cH_p$ the selective compound of hypergraphs $\cH_i$, $i=1,\dots,p$.
Let us note that {\sc Nim}$_{\cH_1\otimes\dots\otimes\cH_p}$ is the conjunctive compound of the games $\Gamma_i=\Nim_{\cH_i}$, $i=1,...,p$, while {\sc Nim}$_{\cH_1\boxplus\dots\boxplus\cH_p}$ is the
selective compound of the same games.

We are now ready to state our main results, which are explicit formulas
for the SG functions of the conjunctive and selective compounds of hypergraph {\sc Nim} games corresponding to JM+ hypergraphs.
We show that JM+ hypergraphs are closed under conjunctive compounds. In fact we prove the following, more general statement.

\begin{theorem}\label{t-main-C}
The conjunctive compound $\cH_1\otimes\dots\otimes\cH_p$, $p\geq 2$ is JM+ if $\cH_i$ is JM+ or $\cH_i=\binom{[2]}{1}$ for all $i=1,\dots,p$. 
\end{theorem}

While the family of JM+ hypergraphs is not closed under selective compounds, we can still describe their SG functions. 

\begin{theorem}
\label{t-main}
Let us consider JM+ hypergraphs $\cH_i\subseteq 2^{V_i}$, $i=1,...,p$, and
their selective compound $\cH=\cH_1\boxplus \dots \boxplus\cH_p$.
For a position $\bx=(\bx^1,...,\bx^p) \in \ZZP^{V_1}\times \cdots \times\ZZP^{V_p}$ of $\Nim_\cH$ let us define
\begin{align*}
\M(\bx) &= m(\bx^1)+\cdots + m(\bx^p),\\
\Y(\bx) &= y_{\cH_1}(\bx^1)+\cdots + y_{\cH_p}(\bx^p), ~\text{and}\\
\H(\bx) &= h_{\cH_1}(\bx^1)+\cdots + h_{\cH_p}(\bx^p).
\end{align*}
Then the SG function of $\Nim_\cH$ is defined by
\[
\G_{\Nim_\cH}(\bx) ~=~ \U(\M(\bx),\Y(\bx),\H(\bx)).
\]
\end{theorem}

Note that if $p=1$, then we get the JM formula.
For both theorems it is an open question if we can replace JM+ hypergraphs by JM ones.

All three considered compounds of
hypergraph {\sc Nim} games provide a hypergraph {\sc Nim} game.
However,  the explicit formula for
the SG function is known only in a  few cases.
In \cite {BGHM15,BBM18} the combined compound
$\Gamma = \Gamma_1 \boxplus (\Gamma_2 \oplus \Gamma_3)$,
where  $\Gamma_i$ is {\sc Nim} with only one pile
for  $i  = 1,2,3$,  was studied.
It is easy to see that this game is the
hypergraph {\sc Nim} on $\cH = \{1,2\},\{1,3\},\{1\},\{2\},\{3\}\}$.
It appears that the SG function $\G_\Gamma(\bx)$  of this game
behaves in a chaotic way when $\bx$  is small and
becomes more regular only for large $\bx$.
Yet, even then no explicit formula for  $\G_\Gamma(\bx)$ is known.

\section{Illustrative examples}
\label{section-example}
Recall that a function $g: \ZZP^V\to \ZZP$ is the {\em SG function} of  {\sc Nim}$_{\cH}$ if and only if the following two conditions hold \cite{Gru39,Spr35,Spr37}.

\begin{itemize}
\item $g(\bx)\neq g(\bx')$  for any move $\bx \to \bx'$;
\item for every integer  $z$ such that
$0 \leq z < g(\bx)$ there exists a move $\bx \to \bx'$ in {\sc Nim}$_\cH$ such that $g(\bx')=z$.
\end{itemize}

Let us recall the SG theorem stating that SG function of the disjunctive compound of impartial games is a function of the SG function values of components (namely, the so-called {\sc Nim}-sum of the SG values \cite{Bou901,Gru39,Spr35,Spr37}). 
Furthermore, in disjunctive compounds, each move to a lower SG value can be realized by moving to a lower SG value in one of the components. 

We will show two examples demonstrating that conjunctive and selective compounds do not have such properties.  
Our first example shows that the SG values of conjunctive and selective compounds are not uniquely defined by the SG values of the components.

\begin{example}
\rm
Let $V_1=\{1,2,3,4\}$ and $V_2=\{5,6,7,8\}$. Define $\cH_1\subseteq 2^{V_1}$ and  $\cH_2\subseteq 2^{V_2}$ by
\[
\cH_1=\{\{1,2\},\{2,3\},\{3,4\},\{4,1\}\} ~\mbox{ and }~\cH_1=\{\{5,6\},\{6,7\},\{7,8\},\{8,5\}\}. 
\]
Then consider positions $\ba^{1}=(0,4,4,0)$ and $\ba^{2}=(0,3,3,0)$ in {\sc Nim}$_{\cH_{1}}$ and {\sc Nim}$_{\cH_{2}}$, respectively.  
For these positions, we have 
\[
m(\ba^{1})=m(\ba^{2})=0, ~y_{\cH_{1}}(\ba^1)=4, \mbox{ and }y_{\cH_{2}}(\ba^2)=3. 
\]
Since both positions are long, we have 
\[
\U(\ba^{1})=h_{\cH_1}(\ba^{1})=4~ \mbox{ and } ~\U(\ba^{2})=h_{\cH_2}(\ba^{2})=3. 
\]
For the pair of positions $\bb^{1}=(0,4,4,0)$ and $\bb^{2}=(4,6,6,4)$, 
we have 
\[
m(\bb^{1})=0, m(\bb^{2})=4, ~y_{\cH_{1}}(\bb^1)=4, \mbox{ and }y_{\cH_{2}}(\bb^2)=2. 
\]
Since $\bb^{1}$ is long and $\bb^{2}$ is short,  
we get 
\[
\U(\bb^{1})=h_{\cH_1}(\bb^{1})=4 ~\mbox{ and }~ \U(\bb^{2})=v(m(\bb^{2}), y_{\cH_2}(\bb^{2}))=3. 
\]
By Theorem \ref{t-main-C}, both hypergraphs are JM+, 
since they are conjunctive compounds of two copies of $\binom{[2]}{1}$.
Consequently, the SG values of $\ba^{i}$ and $\bb^{i}$ are the same for both $i=1,2$.  
 
Let us first consider the conjunctive compound $\cH=\cH_{1} \otimes \cH_2$.
For position $\ba=(\ba^1,\ba^2)$, we have 
\[m(\ba)=\min\{m(\ba^1), m(\ba^2)\}=0  \mbox{ and } 
 y_{\cH}(\ba)=\min\{y_{\cH_1}(\ba^{1}),y_{\cH_2}(\ba^{2})\}=3.
 \]
Since $\ba$ is long,  $\U(\ba)=h_{\cH}(\ba)=3$. 
For position $\bb=(\bb^1,\bb^2)$, we have 
\[m(\bb)=\min\{m(\bb^1), m(\bb^2)\}=0  \mbox{ and } 
 y_{\cH}(\bb)=\min\{y_{\cH_1}(\bb^{1}),y_{\cH_2}(\bb^{2})\}=2.
 \]
Since $\bb$ is also long,  $\U(\bb)=h_{\cH}(\bb)=4$.
By Theorem \ref{t-main-C} $\cH$ is JM+, and thus $\U(\ba)\not=\U(\bb)$ implies that the SG values of $\ba$ and $\bb$ are different.

Let us next consider the selective compound $\cH=\cH_{1} \boxplus \cH_2$. 
Then by applying Theorem \ref{t-main} to  $\cH$,  
we can compute SG values of $\ba=(\ba^1,\ba^2)$ and $\bb=(\bb^1,\bb^2)$ as follows. 

For position $\ba$, we have 
$\M(\ba)=m(\ba^1)+ m(\ba^2)=0$  and 
 $\Y_{\cH}(\ba)=y_{\cH_1}(\ba^{1})+y_{\cH_2}(\ba^{2})=7$.
Thus the SG value $\G(\ba)$ is given by $\U(\ba)=h_{\cH}(\ba)=7$. 
For position $\bb=(\bb^1,\bb^2)$, we have 
$\M(\bb)=m(\bb^1)+ m(\bb^2)=4$ and 
 $\Y(\bb)=y_{\cH_1}(\bb^{1})+y_{\cH_2}(\bb^{2})=6$.
Thus  $\G(\bb)=\U(\bb)=h_{\cH}(\bb)=14$.
This implies that  $\G(\ba)\not=\G(\bb)$. \qed
\end{example}

The next example shows that 
to move to a position with smaller SG value in a selective compound, it may be necessary to increase the SG value in some of the component games. 

\begin{example}
\rm
Let us consider two copies of the hypergraph on $3$ vertices consisting of all proper subsets (Moore's game on $3$ vertices), the positions $\ba^1=\ba^2=(4,4,5)$, 
and the position $\ba=(\ba^1,\ba^2)$ in the selective compound. We have 
\[m(\ba^1)=m(\ba^2)=4, ~y_{\cH_1}(\ba^1)=y_{\cH_2}(\ba^2)=1,\] 
and hence $\U(\ba^1)=\U(\ba^2)=v(4,1)=1$. 
Since both games are JM, we have $\G(\ba^1)=\G(\ba^2)=1$.

In the selective compound, we have
\[\M(\ba)=8, \Y(\ba)=2, \mbox{ and }~ \U(\ba)=v(8,2)=4.
\] 
By Theorem \ref{t-main}, $\U$ is the SG function of the compound game, and hence we must have a move $\ba\to \bb$ such that $\U(\bb)=2$. 
It is easy to argue that for such a position $\bb$ we must have $\Y(\bb)=1$ and $\M(\bb) \equiv 1 \mod 2$. 
The only move (up to the symmetry between the two games) yielding these values is to $\bb^1=(3,3,4)$ in one of the games and  to $\bb^2=(4,4,4)$ in the other one. However the SG value $\G(\bb^1)=\U(\bb^1)=2$ is larger than $\G(\ba^1)=1$.
\qed
\end{example}

\section{Technical Lemmas}
In this section, we present several lemmas which will be used to show our main theorems. 

For positions $\bx,\bx'\in\ZZP^V$ we define
\[
\| \bx'-\bx\|_+ ~=~ \sum_{i\in V\atop x_i'>x_i} x_i'-x_i;
\]
in particular, we have $\| \bx'-\bx\|_+ ~=~  0$ if
$\bx'\leq \bx$, i.e., $x_i'\leq x_i$ for all $i\in V$.

\begin{lemma}\label{l-h-value}
For positions $\bx,\bx'\in\ZZP^V$ we have
\[
h_\cH(\bx)\geq h_\cH(\bx') - \| \bx'-\bx\|_+.
\]
In particular, function $h_\cH$ is monotone with respect to the componentwise relation "$\geq$".
\end{lemma}
\proof
By the definition of the height, it is easy to see that if $\bx\geq \bx'$ then $h_\cH(\bx)\geq h_\cH(\bx')$. Furthermore,
if $\be_j$ is the $j$th unit vector for $j\in V$, then we have $h_\cH(\bx-\be_j) \geq h_\cH(\bx)-1$.
\qed

\medskip

For two integers, $a,b\in \ZZP$ with $a\leq b$,  we denote by $[a,b]$ the set of integers between $a$ and $b$. 
For two positions $\ba,\bb\in\ZZP^V$ with $\ba \leq \bb$, we denote by
\[
[\ba,\bb] ~=~ \{x\in\ZZP^V\mid a_i  \leq x_i \leq b_i, ~ i\in V\}
\]
the set of integer vectors between $\ba$ and $\bb$.

Given a hypergraph $\cH\subseteq 2^V$, a move $\bx\to \bx'$ is called {\em $H$-move} if $x'_i < x_i$ for $i \in H$ and $x_i'=x_i$ for $i \not\in H$. 

\begin{lemma}[Contiguity Lemma]\label{l-h-contiguity}
For a position $\bx\in\ZZP^V$ and an edge $H\in\cH$, let $\bx\to \ba$ and $\bx\to \bb$ denote $H$-moves such that $\ba \leq \bb$.  
Then, all positions $\bc \in [\ba,\bb]$ can be reached by an $H$-move from $\bx$, and we have
\[
\left\{h_\cH(\bc)\mid \bc \in [\ba,\bb]\right\} ~=~ \left[h_\cH(\ba),h_\cH(\bb)\right].
\]
\end{lemma}
\proof
Since both $\bx\to\ba$ and $\bx\to\bb$ are $H$-moves (with the same edge $H\in\cH$), any position $\bc\in [\ba,\bb]$ satisfies  $\bc < \bx$ and $\{i\in V\mid \bc_i < \bx_i\}=H$, proving that $\bx\to\bc$ is an $H$-move. 
Moreover, by the monotonicity of $h_\cH$ in Lemma  \ref{l-h-value}, we have  
$
\left\{h_\cH(\bc)\mid \bc \in [\ba,\bb]\right\} \subseteq  \left[h_\cH(\ba),h_\cH(\bb)\right].
$
To show the converse inclusion, let us define $p=\sum_{i\in V} (b_i-a_i)$, and consider a sequence of positions $\bx^0$,$\bx^1$, \dots, $\bx^p$, such that $\bx^0=\bb$, $\bx^p=\ba$, and for all $j=1,\dots ,p$ $\bx^j$ is obtained from $\bx^{j-1}$ by decreasing one of its components by one unit. Then, again by Lemma \ref{l-h-value} we have $h_\cH(\bx^{j-1})\geq h_\cH(\bx^j)\geq h_\cH(\bx^{j-1})-1$, which proves the converse inclusion.  
\qed

\medskip

Given a hypergraph $\cH\subseteq 2^V$, a move $\bx\to \bx'$ is called  {\em slow} if $x_i-x_i'\leq 1$ for all $i \in V$. 
We denote by $\bx^{s(H)}$ that is obtained from $\bx$ by a slow $H$-move, that is,
\[
x^{s(H)}_i=\begin{cases}
x_i-1&  \text{if } i\in H,\\
x_i & \text{otherwise.}
\end{cases}
\]
\begin{lemma}\label{l-to-long} 
Consider a JM+ hypergraph $\cH$ and a position $\bx\in\ZZP^V$ with $h_\cH(\bx)>0$. 
Then there exists an integer $t$ such that  
\begin{itemize}
\item[(L0)] $t=h_\cH(\bx)$ if $m(\bx)=0$ and $m(\bx)\leq t < h_\cH(\bx)$ if $m(\bx) >0$;
\item[(L1)] for all $t< z<h_\cH(\bx)$ there exists a move $\bx\to \bx'$ such that $0\leq m(\bx')<m(\bx)$, $y_\cH(\bx')\geq y_\cH(\bx)$, and $h_\cH(\bx')=z$;
\item[(L2)] for $z=t < h_{\cH}(\bx)$ there exists a move $\bx\to \bx'$ such that 
$m(\bx')=0$, $y_\cH(\bx')\geq y_\cH(\bx)$, and $h_\cH(\bx')=z$;
\item[(L3)] for all $m(\bx)\leq z< t$ there exists a move $\bx\to \bx'$ such that $m(\bx')=0$ and $h_\cH(\bx')=z$.
\end{itemize}
\end{lemma}

\proof
Let us first consider a position $\bx\in\ZZP^V$ with $m(\bx)>0$. By property (B) there exists a $j\in H\in\cH$ such that $h_\cH(\bx^{s(H)})=h_\cH(\bx)-1$ and $x_j=m(\bx)$. 
Let $\ba^{0}=\bx^{s(H)}$ and define $\bc^{0}$ by   
\[
c^{0}_i~=~\begin{cases} 0&\text{ if } i=j\\x_i-1&\text{ if } i\in H\setminus\{j\}\\
x_i &\text{ if } i\not\in H. \end{cases}
\]
We claim that $t=h_{\cH}(\bc^{(0)})$ has the desired properties. Clearly, this choice satisfies (L0).

For any $\bx' \in [\bc^{0},\ba^{0}]$, there exists a move $\bx \to \bx'$. 
Note that $m(\bx') < m(\bx)$. By this, we have  $\bx'-m(\bx')\be \geq \bx-m(\bx)\be$, 
which implies  that $y_\cH(x') \geq y_\cH(x)$. 
Thus (L2) holds by taking $\bx'=\bc^{0}$, since $m(\bc^{0})=0$. 
Moreover, since $[h_\cH(\bc^{0}),h_\cH(\ba^{0})]=[t,h_\cH(\bx)-1]$, 
Lemma \ref{l-h-contiguity} implies (L1). 

Let us next show (L3). 
By property (A), $\cH_{V\setminus \{1\}}$ has a transversal edge $H'$. 
By property (B), for two edges $H$ and $H'$, 
there exists a chain $H_0\,(=H), H_1, \dots , H_r\,(=H')$. 
Let us then define positions $\ba^{k}$ ($k=1, \dots , r$) by  
\[
a^{k}_i~=~\begin{cases} 0&\text{ if } i\in H_{k-1}\cap H_k\\
x_i-1&\text{ if } i\in H_k\setminus H_{k-1}\\
x_i &\text{ if } i\not\in H_k,\end{cases}
\]
 and $\bb^{k}$ ($k=0, \dots , r$) by 
\[
b^{k}_i~=~\begin{cases} 0&\text{ if } i\in H_k\\x_i &\text{ if } i\not\in H_k. \end{cases}
  \]
Consider the set of positions $I=[\bb^{0}, \bc^{0}] \cup \bigcup_{k=1}^r [\bb^{k}, \ba^{k}]$. 
Note that any position $\bx'$ in $I$ has a move $\bx\to \bx'$ and satisfies $m(\bx')=0$. 
Moreover, for any $i=1, \dots ,r$, we have $\| \bb^{i-1}-\ba^{i}\|_+=1$, implying that 
$h_{\cH}(\ba^{i}) \geq h_{\cH}(\bb^{i-1})-1$. 
By $h_{\cH}(\bb^r) =m$ and Lemma \ref{l-h-contiguity}, it holds that 
\[\{h_{\cH}(\bx')\mid \bx' \in I \}=[m(\bx),t], 
\]
which completes the case of $m(\bx)>0$.

\medskip

Let us finally consider a position $\bx\in\ZZP^V$ such that $m(\bx)=0$. 
We claim that $t=h_{\cH}(\bx)$ satisfies the desired property, i.e., (L3), where (L0), (L1), and   (L2) are automatically satisfied for this $t$.  

Consider $W=\{i\in V\mid x_i>0\}$ and the induced subhypergraph $\cH_W$. By the definition of the height, there exists an edge $H\in\cH_W$ such that $h_\cH(\bx^{s(H)})=h_\cH(\bx)-1$. 
By property (A),  there exists an edge $H'\in\cH_W$ that intersects all other edges of $\cH_W$. 
By property (C),  we have again a chain $H_0\,(=H),H_1,\dots ,H_r\,(=H')$.  Similarly to the above construction, we have a series of $H_k$-moves $k=0,\dots,r$ such that the range of $h_\cH$ values includes all integers $0\leq z < h_\cH(x)$. 
\qed

\begin{lemma}\label{l-to-short}
Assume that $\cH\subseteq 2^V$ satisfies properties (A) and (C). Then for every position $\bx\in\ZZP^V$ with $m(\bx)>0$, 
and pair of integers $(\mu,\eta)\neq(m(\bx),y_\cH(\bx))$ such that 
$0\leq \mu \leq m(\bx)$ and $m(\bx)-\mu \leq \eta \leq y_\cH(\bx)$ there exists a move $\bx\to \bx'$ such that $m(\bx')=\mu$ and $y_\cH(\bx')=\eta$.
\end{lemma}

\proof
Let us consider first the case when $\mu=m(\bx)$. Then (L3) of Lemma \ref{l-to-long} implies the claim when applied to the truncated vector $\bx-m(\bx)\be$.

Let us consider next the case when $\mu < m(\bx)$.
Assume that $x_j=m(\bx)$. Let $H \in \cH$ be an edge with $j\in H$,  and $H' \in \cH$ be a transversal edge of $\cH_{V\setminus \{j\}}$. 
By property (C), we have a chain $H_0\,(=H),H_1,\dots ,H_r\,(=H')$. Define positions $\ba^{k}$ and $\bb^{k}$, $k=0, \dots , r$ by 
\[
a^{k}_i~=~\begin{cases} \mu&\text{ if } i\in H_{k-1}\cap H_k\\
x_i-1&\text{ if } i\in H_k\setminus H_{k-1}\\
x_i &\text{ if } i\not\in H_k\end{cases}
~~~\mbox{ and }~~
b^{k}_i~=~\begin{cases} \mu&\text{ if } i\in H_k\\x_i &\text{ if } i\not\in H_k,  \end{cases}
  \]
where we assume that $H_{-1}=\{j\}$. 
For $k=0, \dots , r$, let $I_k=[\bb^{k}, \ba^{k}]$, and let $I=\bigcup_{k=0}^rI_k$.
We claim that the set of positions $I$ is a certificate of the lemma.  

Since it is clear that any position $\bx' \in I$ has $m(\bx')=\mu$, it is enough to show that 
\begin{equation}
\label{eq-eqa1}
\{y_\cH(\bx')\mid \bx'\in I\}\supseteq [m(\bx)-\mu, y_\cH(\bx)].
\end{equation}    
Note that $y_\cH(\bx')=h_\cH(\bx'-\mu\be)$ for any position $\bx' \in I$. 
Thus for $k=0, \dots , r$, 
 Lemma \ref{l-h-contiguity} implies  that 
\[\{y_{\cH}(\bx')\mid \bx' \in I_k \}=[y_{\cH}(\bb^{k}), y_{\cH}(\ba^{k})].  
\]
Since $y_{\cH}(\ba^{0})\geq y_{\cH}(\bx)$, $y_{\cH}(\bb^{r})=m(\bx)-\mu$, 
and $y_{\cH}(\ba^{k})\geq y_{\cH}(\bb^{k-1})-1$ for $k=1, \dots ,r$, we have 
(\ref{eq-eqa1}), which completes the proof. 
\qed

\begin{lemma}\label{l-selprop}
Assume that $\cH=\cH_1\boxplus\dots\boxplus\cH_p\subseteq 2^V$ is a selective compound of transversal free hypergraphs $\cH_i$, $i=1,\dots,p$. Then $\cH$ itself is transversal free, and for every position $\ba\in\ZZP^V$ and move $\ba\to\bb$ the following relations hold:
\begin{itemize}
\item[(i)] $\H(\ba)>\H(\bb)\geq \M(\ba)\geq\M(\bb)$;
\item[(ii)] $v(\M(\ba),\Y(\ba)) < \M(\ba)$ if and only if $\M(\ba) > \binom{\Y(\ba)+1}{2}$;
\item[(iii)] $\Y(\bb) \geq \M(\ba)-\M(\bb)$.
\end{itemize}
\end{lemma}

\proof
By the definition of the height, it strictly decreases with every move. Moreover, the $m(\ba^i)$ values also can only decrease with a move. To complete the proof of (i), assume that $\ba\to\bb$ is an $H$-move for some $H\in\cH$. By the definition of selective compound, we have $H\cap V_i\in\cH_i\cup\{\emptyset\}$ for all $i=1,\dots,p$. Since all these hypergraphs are transversal free by our assumption, there exists edges $H_i\in\cH_i$ such that $H_i\cap H=\emptyset$ for all $i=1,\dots,p$. Thus even after the $\ba\to\bb$ move we still can make at least $m(\ba^i)$ slow $H_i$-moves from $\bb$. Since these moves for $i=1,...,p$ are all moves in {\sc Nim}$_\cH$, the inequality $\H(\bb)\geq \M(\ba)$ follows. 
The same argument shows also that we can make at least $m(\ba^i)-m(\bb^i)$ slow $H_i$-moves from $\bb$ without decreasing $m(\bb^i)$, for all $i=1,...,p$, proving (iii). 
Finally, (ii) follows by the definition \eqref{e-vmy}. 
\qed

\begin{lemma}\label{l-neq}
Assume that $\cH=\cH_1\boxplus\dots\boxplus\cH_p\subseteq 2^V$ is a selective compound of transversal free hypergraphs $\cH_i\subseteq 2^{V_i}$, $i=1,\dots,p$. Then, for all positions $\bx\in\ZZP^V$ and moves 
$\bx\to \bx'$ of {\sc Nim}$_\cH$ we have 
\[
\U(\M(\bx),\Y(\bx),\H(\bx))\neq \U(\M(\bx'),\Y(\bx'),\H(\bx')).
\]
\end{lemma}

\proof
To prove this statement, we consider four cases, depending on the types of the positions $\bx$ and $\bx'$.  
For simplicity, we use $\U(\bx)$ for $\U(\M(\bx),\Y(\bx),\H(\bx))$. 

If $\M(\bx)\leq \binom{\Y(\bx)+1}{2}$ and $\M(\bx')\leq \binom{\Y(\bx')+1}{2}$, 
then $\U(\bx) = \H(\bx)=h_\cH(\bx)>h_\cH(\bx')= \H(\bx')=\U(\bx')$,
since every move strictly decreases the height by its definition.

If $\M(\bx)\leq \binom{\Y(\bx)+1}{2}$ and $\M(\bx')>  \binom{\Y(\bx')+1}{2}$,  
then we have $\U(\bx)=\H(\bx) > \M(\bx)\geq \M(\bx') > v(\M(\bx'),\Y(\bx'))=\U(\bx')$, proving the claim.
Here the first two inequalities are implied by (i) of Lemma \ref{l-selprop}, while the last inequality follows by (ii) of the same lemma. 

If $\M(\bx)> \binom{\Y(\bx)+1}{2}$ and $\M(\bx')\leq \binom{\Y(\bx')+1}{2}$,   
then we have $\U(\bx')=\H(\bx')\geq \M(\bx)$ by (i) of Lemma \ref{l-selprop}, and
$\M(\bx) > v(\M(\bx),\Y(\bx))=\U(\bx)$ by (ii) of the same lemma.

Finally, if $\M(\bx)> \binom{\Y(\bx)+1}{2}$ and $\M(\bx')> \binom{\Y(\bx')+1}{2}$,  
then by \eqref{e-y} and (i) of Lemma \ref{l-selprop} we have either $\M(\bx)=\M(\bx')$ and $\Y(\bx)>\Y(\bx')$, or $\M(\bx)>\M(\bx')$.

If we have $\Y(\bx)>\Y(\bx')$, then $v(\M(\bx),\Y(\bx))\neq v(\M(\bx'),\Y(\bx'))$ follows by \eqref{e-vmy}.
If $\Y(\bx)=\Y(\bx')$ then by definition \eqref{e-vmy} we could have $v(\M(\bx),\Y(\bx))= v(\M(\bx'),\Y(\bx'))$ if and only if
$\M(\bx')=\M(\bx)-\alpha (\Y(\bx)+1)$ for some positive integer $\alpha$, implying
$\M(\bx')\leq \M(\bx)-\Y(\bx)-1$. From this by (iii) of Lemma \ref{l-selprop} $\Y(\bx')\geq \Y(\bx)+1$ would follow, 
contradicting $\Y(\bx)=\Y(\bx')$. This contradiction shows that we must have $v(\M(\bx),\Y(\bx))\neq v(\M(\bx'),\Y(\bx'))$.
\qed

\section{Proof of Theorem \ref{t-main-C}}

In this section we prove that the family of JM+ games is closed under conjunctive compound. In fact we can show that each of the three properties (A), (B), and (C) are closed under conjunctive compound.

\begin{lemma}\label{l-A-c}
Assume that hypergraphs $\cH^1\subseteq 2^{V_1}$ and $\cH^2\subseteq 2^{V_2}$ satisfy property (A) and $V_1\cap V_2=\emptyset$. Then $\cH=\cH^1\otimes \cH^2$ also satisfies property (A).
\end{lemma}

\proof
Let us note first that for any edge $H=H_1\cup H_2$ of $\cH$ with $H_i\in \cH^i$, $i=1,2$ there exists edges $H_i'\in\cH^i$, $i=1,2$ such that $H_i'\cap H_i=\emptyset$, since $\cH^i$ are transversal free for $i=1,2$.
Thus, $H'=H_1'\cup H_2'\in \cH$ is disjoint from $H$, implying that $\cH$ is transversal free.

Let us consider next a proper subset $S\subsetneq V_1\cup V_2$ such that $\cH_S\neq\emptyset$, and define
$S_i=S\cap V_i$ for $i=1,2$. W.l.o.g., we have $S_1\neq V_1$. Then by the minimal transversal freeness of $\cH^1$ there exists an edge $H_1\in \cH^1_{S_1}$ that intersects all other edges of $\cH^1_{S_1}$. Let us then consider an arbitrary edge $H_2\in\cH^2_{S_2}$. Such an  $H_2$ exists since $\cH_S\neq\emptyset$. Then the set $H=H_1\cup H_2$ is a transversal edge of $\cH_S$.
\qed

\begin{lemma}\label{l-B-c}
Assume that hypergraphs $\cH^1\subseteq 2^{V_1}$ and $\cH^2\subseteq 2^{V_2}$ satisfy property (B) and $V_1\cap V_2=\emptyset$. Then $\cH=\cH^1\otimes \cH^2$ also satisfies property (B).
\end{lemma}

\proof
Let us consider a position $\bx=(\bx^1,\bx^2)\in \ZZP^V$, where $V=V_1\cup V_2$, and note that 
\[
h_{\cH}(\bx) ~=~ \min \{h_{\cH^1}(\bx^1),h_{\cH^2}(\bx^2)\}.
\]
By our assumptions, for both $i=1,2$ we have edges $H_i\in\cH^i$ and $H_i$-moves $\bx^i\to \by^i$ such that 
$h_{\cH^i}(\by^i)=h_{\cH^i}(\bx^i)-1$ and $m(\by^i)<m(\bx^i)$. 
Then with the edge $H=H_1\cup H_2\in\cH$ we can move from $\bx$ to $\by=(\by^1,\by^2)$ and have 
$h_{\cH}(\by)=h_{\cH}(\bx)-1$ and $m(\by)<m(\bx)$. 
\qed

\begin{lemma}\label{l-C-c}
Assume that hypergraphs $\cH^1\subseteq 2^{V_1}$ and $\cH^2\subseteq 2^{V_2}$ satisfy property (C) and $V_1\cap V_2=\emptyset$. Then $\cH=\cH^1\otimes \cH^2$ also satisfies property (C).
\end{lemma}

\proof
Let us consider two edges $H,H'\in\cH$ and denote by $H_i=H\cap V_i$, $H_i'=H'\cap V_i$, and $S_i=H_i\cup H_i'$, for $i=1,2$.
By our assumption there are two chains $F_0,F_1,\dots,F_p$ in $\cH^1_{S_1}$ and $G_0,G_1,\dots,G_r$ in $\cH^2_{S_2}$, such that $F_0=H_1$, $F_p=H_1'$, $G_0=H_2$, and $G_r=H_2'$. Let us then define a chain in $\cH$ such that

\[
F_0\cup G_0, F_1\cup G_0,\dots, F_p\cup G_0, F_p\cup G_1,\dots F_p\cup G_r.
\]
Then we have $H=F_0\cup G_0$, $H'=F_p\cup G_r$, and all the above edges are contained by $H\cup H'$. 
\qed

The following claim can easily be seen.
\begin{lemma}\label{l-D-c}
The conjunctive compound $\cF=\cH\otimes\binom{[k]}{1}$ for some positive integer $k$ satisfies property (C) if $\cH$ satisfies it or $\cH=\binom{[\ell]}{1}$ for some positive integer $\ell$.\qed
\end{lemma}

\begin{theorem}\label{t-2}
The conjunctive compound $\Gamma_1\otimes \Gamma_2$ is a JM+ game, 
whenever $\Gamma_i$ is JM+ or $\Gamma_i=\Nim_2$ for  $i=1,2$. 
\end{theorem}

\proof
If both games are JM+, then the claim follows by Lemmas \ref{l-A-c}, \ref{l-B-c}, and \ref{l-C-c}.
It is also easy to see that the two pile $\Nim_2=\Nim_{\binom{[2]}{1}}$ game satisfies both properties (A) and (B). 
Thus the statement follows by Lemma \ref{l-D-c}. 
\qed

\section{Proof of Theorem \ref{t-main}}

Let us note first that if $p=1$, then the statement is equivalent with saying that minimal transversal free, minimum decreasing hypergraphs that have the chain property are JM. This was already proved in \cite[Theorem 2]{BGHMM17}. Thus we can assume in the sequel that $p\geq 2$,
$\cH=\cH_1\boxplus\dots\boxplus\cH_p\subseteq 2^V$, where  $\cH_i\subseteq 2^{V_i}$, $i=1,\dots,p$, and $V=V_1\cup\dots\cup V_p$. 

To simplify notation we introduce $\U(\bx)=\U(\M(\bx),\Y(\bx),\H(\bx))$ for $\bx\in\ZZP^V$.
To prove the theorem we shall show that the function $\U(\bx)$ satisfies the sufficient conditions for a function to be the SG function, namely that
\begin{itemize}
\item[(D)] for all $\ba\in\ZZP^V$ and moves $\ba\to \bb$ we have $\U(\ba)\neq \U(\bb)$;
\item[(E)] for all $\ba\in\ZZP^V$ and values $0\leq Z < \U(\ba)$ there exists a move $\ba\to\bb$ such that $\U(\bb)=Z$.
\end{itemize}

Property (D) follows by Lemma \ref{l-neq}, since JM+ hypergraphs are all transversal free. To prove property (E), we consider the following two cases:

\begin{itemize}
\item[(E1)] for all positions $\ba\in\ZZP^V$ with $\M(\ba)\leq \binom{\Y(\ba)+1}{2}$ and integers $\M(\ba)\leq Z <\H(\ba)$ there exists a move $\ba\to\bb$ such that $\M(\bb)\leq \binom{\Y(\bb)+1}{2}$ and $\H(\bb)=Z$;
\item[(E2)] for all positions $\ba\in\ZZP^V$ and values $0\leq Z < \min(\M(\ba),v(\M(\ba),\Y(\ba)))$ there exists a move $\ba\to\bb$ such that 
$\M(\bb)> \binom{\Y(\bb)+1}{2}$ and $v(\M(\bb),\Y(\bb))=Z$.
\end{itemize}
It is easy to see that properties (E1) and (E2) will imply property (E) by \eqref{e-long} and \eqref{e-short}.

To prove (E1) let us consider a position $\ba=(\ba^1,\dots,\ba^p)\in\ZZP^V$ with $\M(\ba)\leq \binom{\Y(\ba)+1}{2}$. By Lemma \ref{l-to-long} there exist thresholds $m(\ba^i)\leq t_i \leq h_{\cH_i}(\ba^i)$, $i=1,...,p$ satisfying the claims of the lemma. 
Let us set $T=t_1+\dots +t_p$. 

For an integer $T\leq Z <\H(\ba)$ let us choose integers $t_i\leq z_i \leq h_{\cH_i}(\ba^i)$ such that $Z=z_1+\dots +z_p$. 
Let us note that by the above definitions, we must have $m(\ba^i)=0$ whenever $t_i=z_i=h_{\cH_i}(\ba^i)$. 
Let us define $Q=\{i\in [p]\mid z_i<h_{\cH_i}(\ba^i)\}$ and note that $Q\neq\emptyset$ and $m(\ba^i)>0$ for all $i\in Q$. 
Thus, by (L1) and (L2) of Lemma \ref{l-to-long} for every $i\in Q$ there exists a move $\ba^i\to \bb^i$ such that $m(\bb^i)\leq m(\ba^i)$, $y_{\cH_i}(\bb^i)\geq y_{\cH_i}(\ba^i)$, and $h_{\cH_i}(\bb^i)=z_i$. 
We define $\bb^i=\ba^i$ for $i\in [p]\setminus Q$, and set
$\bb=(\bb^1,\dots,\bb^p)$. In this way we get that $\ba\to \bb$ is a move in {\sc Nim}$_\cH$, satisfying $\M(\bb)\leq \M(\ba)$, $\Y(\bb)\geq \Y(\ba)$, and $\H(\bb)=Z$. Thus, by our assumptions, it follows that $\M(\bb)\leq \M(\ba)\leq \binom{\Y(\ba)+1}{2}\leq \binom{\Y(\bb)+1}{2}$. 

For an integer $\M(\ba)\leq Z< T$ let us choose integers $m(\ba^i)\leq z_i\leq t_i$ such that $Z=z_1+\dots +z_p$, and  
define $Q=\{i\in [p]\mid z_i<h_{\cH_i}(\ba^i)\}$ as above. We have $Q\neq\emptyset$, because $Z<T$. Furthermore, for $i\in[p]\setminus Q$ we have $t_i=h_{\cH_i}(\ba^i)$, implying $m(\ba^i)=0$ by (L0) of Lemma \ref{l-to-long}.
By (L2) and (L3) of Lemma \ref{l-to-long} for every $i\in Q$ there exists a move $\ba^i\to \bb^i$ such that $m(\bb^i)=0$ and $h_{\cH_i}(\bb^i)=z_i$. 
We define $\bb^i=\ba^i$ for $i\in [p]\setminus Q$, and set
$\bb=(\bb^1,\dots,\bb^p)$.
Then $\ba\to \bb$ is a move in {\sc Nim}$_\cH$ satisfying $\M(\bb)=0$ and $\H(\bb)=Z$. Thus, it trivially follows that $0=\M(\bb)\leq \binom{\Y(\bb)+1}{2}$. This completes the proof of property (E1). 

\medskip

For the proof of property (E2) we need to make a few more observations. Let us note first that by \eqref{e-vmy} for a fixed integer $y\in\ZZP$ we have
\[
U(y) = \{ v(m,y)\mid m\in\ZZP\} ~=~ \left[\binom{y+1}{2},\binom{y+1}{2}+y\right]
\]
and that the sets $U(y)$, $y\in \ZZP$ partition the set of nonnegative integers. Consequently, for every integer $z\in\ZZP$ there exists a unique integer $y\in\ZZP$ such that $z\in U(y)$. We denote this unique integer as $y=\eta(z)$. Let us also note that for every integer $z$ we have
\begin{equation}\label{e-z+1}
z ~=~ v(z+1,\eta(z)) ~~~\text{ and }~~~ z+1 ~>~ \binom{\eta(z)+1}{2}.
\end{equation}

Let us now consider a position $\ba\in\ZZP^V$ and a value $0\leq Z < \min(\M(\ba),v(\M(\ba),\Y(\ba)))$, as in (E2), and 
choose a largest integer $\alpha\geq 0$ such that $M=Z+1+\alpha\dot(\eta(Z)+1) \leq \M(\ba)$ and set $Y=\eta(Z)$. Note that we have 
\begin{equation}\label{e-M}
0\leq \M(\ba)-M\leq \eta(Z).
\end{equation}
We construct a position $\bb\in\ZZP^V$ such that $\ba\to\bb$ is a move in {\sc Nim}$_\cH$, $\M(\bb)=M$, and $\Y(\bb)=Y$. By \eqref{e-z+1} and the fact that $\alpha\geq 0$ this construction verifies property (E2), and completes our proof of the theorem. 

\smallskip

To see the construction, let us start first by noting that $M\leq \M(\ba)=\sum_{i=1}^p m(\ba^i)$ by our choices above. Thus there exists integer values $\mu_i$, $i=1,\dots,p$ such that 
\[
\sum_{i=1}^p \mu_i ~=~ M, \text{ and } 0\leq \mu_i\leq m(\ba^i) ~\forall~ i=1,\dots p.
\]
Let us next observe that we have 
\[
\sum_{i=1}^p (m(\ba^i)-\mu_i) ~=~ \M(\ba)-M ~\leq~ Y ~\leq~ \Y(\ba) ~=~ \sum_{i=1}^p y_{\cH_i}(\ba^i)
\]
by \eqref{e-M} and by our definitions \eqref{e-m} and \eqref{e-y}. 
Thus there exists integers $\eta_i$, $i=1,\dots,p$ such that 
\[
m(\ba^i)-\mu_i ~\leq~ \eta_i ~\leq~ y_{\cH_i}(\ba^i)~~~\forall i=1,\dots,p,
\]
and $Y=\sum_{i=1}^p \eta_i$. 
Let us define $Q=\{i\in [p]\mid (\mu_i,\eta_i)\neq (m(\ba^i),y_{\cH_i}(\ba^i))\}$. Since by (E2) we have $v(M,Y)=Z<v(m(\ba),\Y(\ba))$ implying $(M,Y)\neq (m(\ba),\Y(\ba))$, and thus $Q\neq\emptyset$.

Now we can apply Lemma \ref{l-to-short} for each of the JM+ games $\cH_i$, $i\in Q$, and derive the existence of moves $\ba^i\to\bb^i$ in {\sc Nim}$_{\cH_i}$ such that 
$m(\bb^i)=\mu_i$ and $y_{\cH_i}(\bb^i)=\eta_i$.
Defining $\bb^i=\ba^i$ for $i\in [p]\setminus Q$ and setting
$\bb=(\bb^1,\dots,\bb^p)$ completes our construction and the proof of the theorem.
\qed

\section{Combinatorial Properties of JM+ Hypergraphs}
\label{s6}

In this section we show that no two of the three properties (A), (B), and (C) imply JM, and in particular the third property. This provides further justification of why JM+ is an interesting subfamily of JM. 

For a hypergraph $\cH\subseteq 2^V$ we denote by $\min \cH$ the family of inclusionwise minimal edges of $\cH$. For a subhypergraph $\cF\subseteq \cH$ we denote by $V(\cF)=\bigcup_{F\in \cF} F$ the set of vertices that it covers.

Let us strengthen property (B) by the following combinatorial property (B*). This is because it will be technically easier to verify in some of the next examples and families than (B) itself, since property (B) involves "$\forall \bx \in\ZZP^V$".

\begin{itemize}
\item[(B*)] For every subhypergraph $\cF\subseteq \min\cH$ such that $V(\cF)\neq V$ there exist edges $F\in\cF$ and $H\in\cH$ such that $H\cap V(\cF)\subseteq F$ and $H\setminus F\neq \emptyset$.
\end{itemize}

A move $\bx \to \bx'$ is called {\em height move} if $h_{\cH}(\bx')=h_{\cH}(\bx)-1$. 

\begin{lemma}\label{l-B->B*}
If a hypergraph $\cH$ satisfies property (B*) then it also satisfies (B). 
\end{lemma}

\proof
Let us consider a position $\bx\in\ZZP^V$ with $m(\bx)>0$ and define $\cF(\bx)\subseteq \cH$ to be the subhypergraph of those edges $H\in\cH$ for which there exists a $\bx\to\bx'$ $H$-move such that $h_\cH(\bx')=h_\cH(\bx)-1$. If $V(\cF(\bx))=V$ for all positions $\bx\in\ZZP^V$, then property (B) holds. Otherwise there exists a position $\bx$ with $m(\bx)>0$ such that $V(\cF(\bx))\subsetneq V$. Note that w.l.o.g. we can assume that $\cF\subseteq \min\cH$. Thus, by property (B*) we have edges $F\in\cF(\bx)$ and $H\in\cH$ such that $H\cap V(\cF(\bx))\subseteq F$ and $H\setminus V(\cF(\bx))\neq \emptyset$. Clearly, we can assume that $H\in\min\cH$. Then there exists a sequence of height moves that involves $F$ by the definition of $\cF(\bx)$. In this sequence let us replace one $F$-move by an $H$-move. This way we get another height sequence, and that contradicts the fact that $H\not\in\cF(\bx)$. This contradiction proves our claim.
\qed

Let us add that the inverse implication is not true. The following small example shows this. Consider $V=\{1,2,3,4,5\}$, $H_1=\{1,2\}$, $H_2=\{2,3\}$, $H_3=\{3,4\}$, $H_4=\{1,4,5\}$, and $\cH=\{H_1,H_2,H_3,H_4\}$. For the subhypergraph $\cF=\{H_1,H_2,H_3\}$ property (B*) fails to hold. It is however not difficult to see that $\cH$ satisfies property (B). 

For a subset $S\subseteq [n]$ we denote by $\binom{[n]}{S}$ the hypergraph consisting of all edges $H\subseteq V$ such that $|H|\in S$. If $S=\{i\}$, then we simply write $\binom{[n]}{i}$. These hypergraphs are called \emph{symmetric}. 
We say that $S\subseteq [n]$ has a gap, if there are integers $0<i<j<k\leq n$ such that $i,k\in S$ and $j\not\in S$.

\begin{remark}
It is easy to see that symmetric hypergraphs satisfy property (B*). Furthermore, if $S$ has a gap, then $\binom{[n]}{S}$ does not satisfy property (C). We also recall from \cite{BGHMM18} that symmetric JM hypergraphs have a simple arithmetic characterization. 
\end{remark}

\begin{figure}[h]
\begin{center}
\begin{tikzpicture}[scale=1]

\filldraw[color=blue,opacity=0.5] (4,2) -- (12,2) -- (12,8) -- (4,8) -- cycle;
\filldraw[color=red,opacity=0.5] (2,4) -- (8,4) -- (8,10) -- (2,10) -- cycle;
\filldraw[color=yellow,opacity=0.5] (0,0) -- (7,0) -- (7,6) -- (0,6) -- cycle;
\filldraw[color=cyan,opacity=0.5] (4,3) -- (10,3) -- (10,7) -- (4,7) -- cycle;

\node at (5.5,5) [shape=circle,draw=none] (JM+) {\tiny $\mathbf{JM+}$};
\node at (1,5.5) [shape=circle,draw=none] (B) {\tiny $(\mathbf{B})$};
\node at (11,7.5) [shape=circle,draw=none] (A) {\tiny $(\mathbf{A})$};
\node at (3,9.5) [shape=circle,draw=none] (C) {\tiny $(\mathbf{C})$};
\node at (9,6.5) [shape=circle,draw=none] (JM) {\tiny $(\mathbf{JM})$};

\node at (11,5) [shape=circle,draw=none] (P1) {\tiny $(\mathbf{P1})$};
\node at (2,2) [shape=circle,draw=none] (P2) {\tiny $(\mathbf{P2})$};
\node at (5,9) [shape=circle,draw=none] (P3) {\tiny $(\mathbf{P3})$};

\node at (5.5,2.5) [shape=circle,draw=none] (P4) {\tiny $(\mathbf{P4})$};
\node at (3,5) [shape=circle,draw=none] (P5) {\tiny $(\mathbf{P5})$};
\node at (6,7.5) [shape=circle,draw=none] (P6) {\tiny $(\mathbf{P6})$};

\node at (5.5,3.5) [shape=circle,draw=none] (P7) {\tiny $(\mathbf{P7})$};
\node at (6,6.5) [shape=circle,draw=none] (P8) {\tiny $({\color{yellow}\mathbf{P8}})$};
\node at (9,5) [shape=circle,draw=none] (P9) {\tiny $({\color{yellow}\mathbf{P9}})$};

\end{tikzpicture}
\end{center}
\caption{\label{figp1p9} The $10$ regions defined by properties (A), (B), (C), and (JM).}
\end{figure}
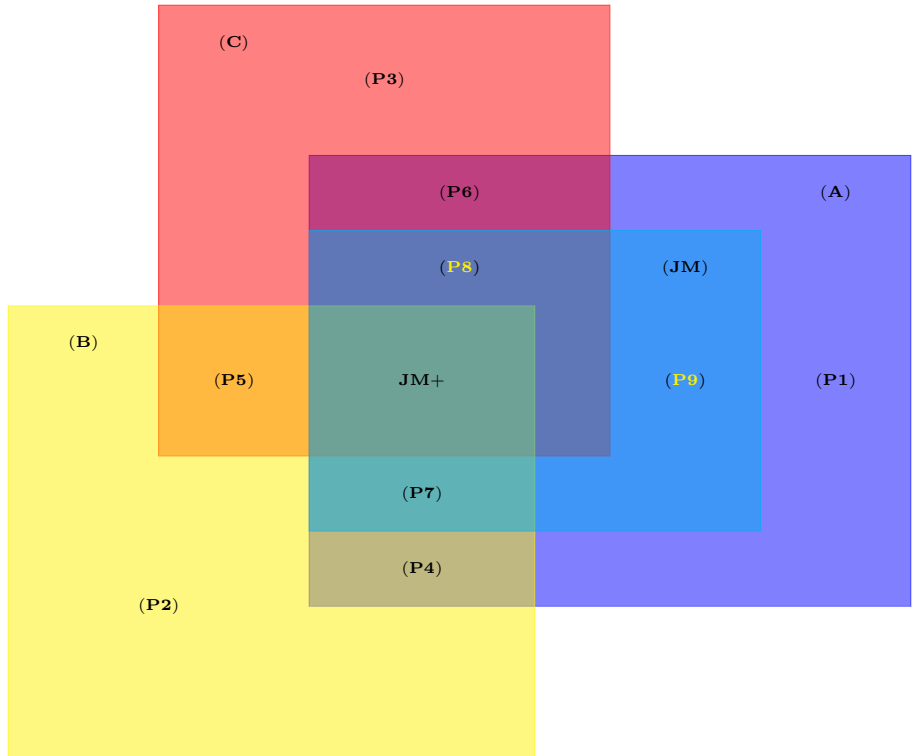

Due to our results, properties (A), (B), (C), and JM define 10 possible regions (see Figure \ref{figp1p9}), one of which are the JM+ hypergraphs. We show below that 7 of the remaining regions are nonempty. The status of the last two is open. 

\begin{theorem}
The following 7 statements are all satisfied. 
\begin{itemize}
\item[(P1)] Property (A) implies none of JM, (B), and (C).
\item[(P2)] Property (B) implies none of (A) and (C).
\item[(P3)] Property (C) implies none of (A) and (B).
\item[(P4)] Properties (A) and (B) do not imply JM, and hence (C).
\item[(P5)] Properties (B) and (C) do not imply JM, and hence (A).
\item[(P6)] Properties (A) and (C) do not imply JM, and hence (B).
\item[(P7)] Properties JM and (B) do not imply (C).
\end{itemize}
\end{theorem}

\proof
For (P1) we consider the "cube" defined on the 8 vertices of a 3 dimensional unit cube, in which the edges are formed by the 4 vertices of 2 dimensional faces. It was shown in \cite[Section 7]{BGHMM17} that this hypergraph satisfies property (A), but none of the others. 

For (P2) we consider e.g., $\binom{[4]}{\{1,3\}}$. As we noted above this satisfies (B*), and hence (B), but it does not satisfy (C). It does ot satisfy (A) either, since if $S$ is a subset of size $2$, then the induced subhypergraph does not have a transversal edge. 

For (P3) we consider the following hypergraph on 10 vertices: define $T_1=\{v_1,v_2,v_3\}$, $T_2=\{v_4, v_5, v_6\}$, $T_3=\{v_7,v_8,v_9\}$, and $V=\{v_0,v_1,\dots, v_9\}$. Consider 
\[
\cH=\{T_1,T_2,T_3\}\cup 
\left\{
H\subseteq V
\left| 
\begin{array}{l}
|H|=4 \text{ and } v_0\not\in H, \text{ or}\\
|H|=5 \text{ and } v_0\in H
\end{array}
\right. \right\}
\]
This hypergraph does not satisfy (A), since if we delete $v_0$ then it still does not have a transversal edge. 
It does not satisfy property (B) either, because for position $\bx=(1,2,2,2,2,2,2,2,2,2)$ the only height moves are with sets $T_i$, $i=1,2,3$. Finally, it is easy to verify that it satisfies property (C). 

For claims (P4), (P5), and (P7) it is enough to consider the symmetric hypergraphs 
$\binom{[2]}{1}$, $\binom{[5]}{2}$, and $\binom{[6]}{\{2,4\}}$, respectively.

Finally, for (P6) we consider a 10-vertex hypergraph similar to the one considered for (P3).
\[
\cH=\{T_1,T_2,T_3\}\cup 
\left\{
H\subseteq V
\left| 
\begin{array}{l}
|H|=4,~ v_0\not\in H, \text{ and } H\cap \{v_1,v_4,v_7\}\neq\emptyset, \text{ or}\\
|H|=5 
\end{array}
\right. \right\}
\]
To see property (A), note first that $\cH$ has no transversal edge. 
Furthermore, if $S\subseteq V$ has $|S|\leq 5$, then any two edges of $\cH_S$ intersect. If $|S|\geq 6$, then any edge $H\in\cH_S$ such that $|H|=5$ and $H \supseteq \{v_1,v_4,v_7\} \cap S$ is a transversal of $\cH_S$. 

For property (C), consider two edges $H, H' \in \cH$. 
If $|H \cup H'| \leq 4$, then it is easy to see that a chain exists from $H$ to $H'$. 
On the other hand, if $|H \cup H'| \geq 5$, choose an edge $H'' \subseteq H \cup H'$ of size $5$. 
Note that we can reach $H'$ from $H$ by a chain through $H''$. 
Thus $\cH$ satisfies (C). 

To show that $\cH$ is not JM, let us consider the position $\bx=(6,7, \dots ,7)$. 
Note that $m(\bx)=6$, $y_\cH(\bx)=3$, and $h_\cH(\bx)=21$. 
Since $\binom{y_\cH(\bx)+1}{2}=m(\bx)$, we have $\U(\bx)= h_\cH(\bx)=21$. 
We will show that there exists no move from $\bx\to \bx'$ such that $\U(\bx')=20$. 
Assume that such a move exists. 
Since $m(\bx')\leq m(\bx)=6$, position $\bx'$ is  long with  $h_\cH(\bx')=20$. 
This implies that $\bx\to \bx'$ is a height move. 
For this position $\bx$, the only height moves are with $T_i$, $i=1,2,3$. 
Consequently, $y_\cH(\bx')< 3$ and $m(\bx')=6$, implying that $\bx'$ is a short position. 
Thus $\bx'$ is both short and long,  a contradiction. 
\qed

For the remaining two statements, we do not know any examples:
\begin{itemize}
\item[(P8)] Properties JM and (C) do not imply (B).
\item[(P9)] Property JM  implies neither (B) nor (C).
\end{itemize}

In fact, we do not know if property JM implies (B) or not. 

\section*{Acknowledgements}
The authors also thank Rutgers University and RUTCOR for the support
to meet and collaborate.
The second author was partially funded by
the Russian Academic Excellence Project '5-100'.


\begin{thebibliography}{99}

\bibitem{Alb07}
M.H. Albert, R.J. Nowakowski, D. Wolfe,
Lessons in play: An introduction to combinatorial game theory,
second ed., A K Peters Ltd., Wellesley, MA, 2007.


\bibitem{Ber62}
C. Berge, The theory of graphs, London, 1962.

\bibitem{BCG01-04}
E.R. Berlekamp, J.H. Conway, and R.K. Guy,
Winning ways for your mathematical plays,
vol.1-4, second edition, A.K. Peters, Natick, MA, 2001 - 2004.

\bibitem{BBM18}
C. Beideman, M. Bowen, N.A. M\"uyesser,
The Sprague-Grundy function for some selective compound games,
arXiv:1802.08700v1  23 Feb. 2018;  https://arxiv.org/abs/1802.08700



\bibitem{BGHM15}
E. Boros, V. Gurvich, Nhan Bao Ho, K. Makino,
On the Sprague-Grundy Function of Tetris Extensions of Proper {\sc Nim},
RUTCOR Research Report, RRR-1-2015, Rutgers University,
available at http://arxiv.org/abs/1504.06926.

\bibitem{BGHMM15}
E. Boros, V. Gurvich, N.B. Ho, K. Makino, and P. Mursic,
On the Sprague–Grundy function of Exact $k$-Nim,
Discrete Appl. Math., 239 (2018) 1--14.

\bibitem{BGHMM16}
E. Boros, V. Gurvich, N.B. Ho, K. Makino, and P. Mursic,
Hypergraph Combinations of Impartial Games with
Decreasing Sprague-Grundy Functions,
available at   http://arxiv.org/abs/1701.02819

\bibitem{BGHMM17}
E. Boros, V. Gurvich, Nhan Bao Ho, K. Makino, and P. Mursic,
On the Sprague-Grundy function of matroids and related hypergraphs,
available at  https://arxiv.org/abs/1804.03692

\bibitem{BGHMM18}
E. Boros, V. Gurvich, Nhan Bao Ho, K. Makino, and P. Mursic,
Sprague-Grundy function of symmetric hypergraphs,
\emph{J. Combinatorial Theory Ser A}, 165 (2019) 176-186.

\bibitem{Bou901}
C.L. Bouton, Nim, a game with a complete mathematical theory,
Ann. of Math., 2-nd Ser. 3 (1901-1902) 35-39.

\bibitem{Con76}
J.H. Conway,  On numbers and games,
Acad. Press, London, NY, San Francisco, 1976.

\bibitem{ES96}
R. Ehrenborg and E. Staingrimsson,
Playing Nim on Simplicial Complex,
The Electronic J. Cominatorics 3 (1996) \#R9.


\bibitem{Gru39}
P.M. Grundy, Mathematics of games, Eureka 2 (1939) 6-8.

\bibitem{GS56}
P.M. Grundy and C.A.B. Smith,
Disjunctive games with the last player losing,
Proc. Cambridge Philos. Soc., 52 (1956) 527-523.

\bibitem{Har69}
F. Harary, Graph Theory, Addison-Wesley, 1969.

\bibitem{JM80}
T.A. Jenkyns and J.P. Mayberry,
Int. J. of Game Theory 9 (1) (1980) 51--63,
The skeletion of an impartial game
and the Nim-Function of Moore's Nim$_k$.


\bibitem{Moo910}
E.H. Moore,
A generalization of the game called {\sc Nim},
Annals of Math., Second Series, 11:3 (1910) 93--94.

\bibitem{Sie13}
A. N. Siegel, Combinatorial Game Theory,
Graduate Studies in Mathematics 146, AMS, 2013.




\bibitem{Smi66}
C.A.B. Smith, Graphs and composite games,
J. of Combinatorial theory 1  (1966) 51--81.

\bibitem{Spr35}
R. Sprague, \"Uber mathematische Kampfspiele,
Tohoku Math. J. 41 (1935-36) 438-444.

\bibitem{Spr37}
R. Sprague, \"Uber zwei abarten von nim, Tohoku Math. J. 43 (1937) 351--354.

\bibitem{Tut54}
W. Tutte, A short proof of the factor theorem for finite graphs,
Canadian Journal of Mathematics, 6 (1954) 347--352.


\end{thebibliography}
\end{document}